\documentstyle[amstex,amscd]{article}

\newtheorem{thm}{Theorem}[section]
\newtheorem{lm}[thm]{Lemma}
\newtheorem{cor}[thm]{Corollary}
\newtheorem{prop}[thm]{Proposition}
\newtheorem{dfn}[thm]{Definition}
\newtheorem{con}[thm]{Conjecture}

\newtheorem{remark}[thm]{Remark}

\newcommand{\bt}{\begin{thm}}
\newcommand{\et}{\end{thm}}
\newcommand{\bl}{\begin{lm}}
\newcommand{\el}{\end{lm}}
\newcommand{\bc}{\begin{cor}}
\newcommand{\ec}{\end{cor}}
\newcommand{\bp}{\begin{prop}}
\newcommand{\ep}{\end{prop}}
\newcommand{\bd}{\begin{dfn}}
\newcommand{\ed}{\end{dfn}}
\newcommand{\bcon}{\begin{con}}
\newcommand{\econ}{\end{con}}
\newcommand{\brem}{\begin{remark}}
\newcommand{\erem}{\end{remark}}

\newcommand{\pf}{\noindent {\bf Proof:}\hspace{2.5mm}}
\newcommand{\qed}{~\hspace{2.5mm}~\rule{2.5mm}{2.5mm}}

\newcommand{\nR}{\text{\normalshape R}}  
\newcommand{\nS}{\text{\normalshape S}}

\newcommand{\gin}{\text{\normalshape gin}}
\newcommand{\ini}{\text{\normalshape in}}
\newcommand{\isat}{{\it sat}}

\newcommand{\grl}{>_{\text{\normalshape rlex}}}

\newcommand{\kk}{{\Bbb K}}

{\end{array}}

\begin{document}

\title{Monomial invariants in codimension two
\footnote{
\newline
The authors are members of CNR--GNSAGA (Italy). 
\newline
During the preparation of this paper the authors were partially supported by
National Research Project ``Geometria sulle variet\`a algebriche COFIN 2002'' 
of MIUR--Italy. }}
\author{Alberto Alzati
\\Alfonso Tortora
\\Dipartimento di Matematica, Universit\`{a} di Milano}
\date{
}
\maketitle


\begin{abstract}
We define the {\it monomial invariants} of a projective variety $Z$; 
they are invariants coming from the generic initial ideal of $Z$.\\ 
Using this notion, we generalize a result of Cook \cite{c}:\\ 
If $Z$ is an integral variety of codimension two, satisfying the
additional hypothesis $s_Z=s_\Gamma,$ then its monomial invariants are 
connected. 
\end{abstract}


\section*{Introduction}
The generic initial ideal of a projective variety $Z,\;\gin(I_Z),$ preserves 
some informations about $Z,$ in particular it has the same Hilbert function; 
on the other hand, $\gin(I_Z)$ is a combinatorial object, which can be 
``finitely'' described in terms of its monomial generators; thus any limitation
on the ``shape'' of $\gin(I_Z),$ i.e. any relation among its generators, 
translates into a limitation on the possible Hilbert functions of projective 
varieties.\\
Generic initial ideals are particularly well suited to study codimension two
varieties; in this situation $\gin(I_Z)$ seems to have just the right amount of
information. Since we shall mostly deal with sections of $Z,$ we use the 
reverse lexicographic order on monomials; also, $\gin(I_Z)$ being a saturated 
ideal, its monomial generators do not contain the last variable $x_n$---here
$Z$ is a (nondegenerate) subvariety of ${\Bbb P}^n.$\\
The first instance of such a variety is a set $\Gamma$ of points in 
${\Bbb P}^2,$ for which $\gin(I_\Gamma)\subseteq\kk[x_0,x_1,x_2]$ is minimally 
generated by monomials of type $x_0^s, x_0^{s-1}x_1^{\lambda_{s-1}},..., 
x_0x_1^{\lambda_1}, x_1^{\lambda_0};$ a classical result of Gruson and Peskine
\cite{gp} implies that, if the points of $\Gamma$ are in general position, then
$\lambda_{i+1}+1\geq\lambda_i\geq\lambda_{i+1}+2,$ for all $i=0, 1,..., s-2;$ 
we describe this situation by saying that the invariants $\lambda_i$ are 
connected, i.e. in passing from a generator $x_0^ix_1^{\lambda_i}$ to the next
$x_0^{i+1}x_1^{\lambda_{i+1}},$ the invariant $\lambda_i$ can ``jump'' downward
to $\lambda_{i+1}$ by one or two steps, but no more (the fact that it jumps no 
less than that is a consequence of another property of $\gin,$ its 
Borel--fixedness).\\
Next case is a space curve $C\subseteq{\Bbb P}^3.$ Now $\gin(I_C)$ can be 
thought of as a stack of slices in the $x_2$ direction, the slice at level 
$p\geq0$ being the set of all monomials $x_0^\alpha x_1^\beta x_2^p$ belonging 
to $\gin(I_C).$ Since the slice at level $p\gg0$ represents $\gin(I_\Gamma),$ 
where $\Gamma$ is the general plane section of $C,$ then the invariants of such
a slice are connected; the question is what happens when $p$ is small. The 
main result of a paper of Cook \cite{c} is that the invariants of $\gin(I_C)$ 
are connected at all levels $p\geq0.$ Unfortunately, in the proof of this 
result there are some gaps. 
In order to fix those gaps, Decker and Schreyer made
the hypothesis that $s_C=s_\Gamma,$  where $s_C$ (or $s_\Gamma$) is the
minimal degrees of a polynomial vanishing on $C$ (or $\Gamma$); this hypothesis
seems to be unavoidable, but the proof given in \cite{ds} is still incomplete.
Finally Amasaki \cite{a} gave a complete proof of the connectedness, under the 
same hypothesis, using  mostly algebraic techniques.\\
The present paper translates Amasaki's ideas in a more geometric language,
generalizing the result about the connectedness of invariants to higher 
dimensional varieties, in the following sense.\\
Let $Z\subseteq{\Bbb P}^n$ be an integral (i.e. reduced and irreducible) 
nondegenerate variety of codimension two and suppose that $s_Z=s_\Gamma,$ 
then the slice of $\gin(I_Z)$ at 
level $p_j$ with respect to the variable $x_j,\; j=2,...,n-1,$ has connected 
invariants, for all $p_j\geq0.$\\ 
{\bf Acknowledgments.} The authors had fruitful discussions--either 
in person or via e-mail--on the topics of this paper with many people, among 
them Matsumi Amasaki, Iustin Coanda, Wolfram Decker, Mark Green, 
Francesco Russo, Enrico Sbarra, Enrico Schlesinger, Frank Schreyer; 
their help is gratefully acknowleged.
 

\section{Background}
Let $\kk$ be the base field and assume that it be algebraically closed.\\
Let $W$ be a $\kk$-vector space of dimension $n+1$ and denote by 
$\kk\langle x_0,...,x_n\rangle$ the $(n+1)$-dimensional $\kk$-vector space 
generated by $x_0,...,x_n.$\\
Our enviroment is the projective space ${\Bbb P}^n={\Bbb P}(W^*),$
whose ring of polynomials is the symmetric algebra 
$\nS W=\oplus_{d\geq0}\nS^dW.$ Of course, starting with
$\kk\langle x_0,...,x_n\rangle,$ we get ${\Bbb P}^n(\kk)$ and 
$\kk[x_0,...,x_n]$ respectively. We use italic $I, J,...$ for ideals in 
$\nS W$
and gothic $\frak{i}, \frak{j},...$ for ideals in $\kk[x_0,...,x_n].$\\
Another standard piece of notation is the multiindex one: 
${\bf x}^{\bf p}=x_0^{p_0}...x_n^{p_n},$ with
$|{\bf p}|=\deg{\bf x}^{\bf p}=p_0+...+p_n.$ 
\bd
A {\em coordinate system} (or simply {\em coordinates}) on $W$ is an 
isomorphism of $\kk$-vector spaces $W\rightarrow\kk\langle x_0,...,x_n\rangle.$
We denote by $\nR(W)$ the set of all coordinate systems on $W.$ 
\ed
Clearly, a coordinate system on $W$ gives coordinates in $\nS W,$ i.e. an 
isomorphim of $\kk$-algebras $\nS W\simeq\kk[x_0,...,x_n].$ Since 
$W\simeq\kk^{n+1},$ there is a bijection $\nR(W)\leftrightarrow GL(n+1,\kk),$
thus $\nR(W)$ has a natural structure of algebraic variety, with its attendant 
Zariski topology.
\bd
(i) We consider an ordering of the set of monomials of $\kk[x_0,...,x_n],$ the
{\em reverse lexicographic order}, denoted by $\grl,$ and defined as follows:\\
${\bf x}^{\bf p}\grl{\bf x}^{\bf q}$ if either 
$\deg{\bf x}^{\bf p}<\deg{\bf x}^{\bf q},$ or 
$\deg{\bf x}^{\bf p}=\deg{\bf x}^{\bf q}$ and 
$p_n=q_n,\ldots,p_{k+1}=q_{k+1} ,p_k<q_k,$ for some $k.$\\
(ii) If $f\in\kk[x_0,...,x_n],$ the {\em initial monomial} of $f, \ini(f),$ 
is the greatest monomial, with respect to $\grl,$ that appears in $f.$\\
(iii) Let ${\frak i}$ be a homogeneous ideal in $\kk[x_0,...,x_n],$ we define 
the {\em initial ideal} of ${\frak i},$ denoted by $\ini({\frak i}),$ as the 
monomial ideal generated by all $\ini(f), f\in{\frak i}.$\\
(iv) The {\em $j$-th elementary move}, $j=1,...,n$ is defined on monomials by\\
$e_j({\bf x}^{\bf p})=
x_0^{p_0}...x_{j-1}^{p_{j-1}+1}x_j^{p_j-1}...x_n^{p_n},$ if $p_j>0,$
otherwise $e_j({\bf x}^{\bf p})=0.$\\
(v) A monomial ideal ${\frak b}$ is {\em Borel-fixed} if, for all monomials 
${\bf x}^{\bf p}\in{\frak b}$ and 
all elementary moves $e_j, \; e_j({\bf x}^{\bf p})\in{\frak b}.$  
\ed
If $I$ is a homogeneous ideal in $\nS W$ and $g\in\nR(W),$ then $gI$ 
is an ideal of $\kk[x_0,...,x_n],$ so it makes sense to consider 
$\ini_g(I):=\ini(gI).$ 
Although $\ini_g(I)$ depends, by its very definition,  on the choice of $g,$ it
turns out that, for general coordinates $g\in\nR(W),\;\ini_g(I)$ stays 
constant, i.e. does not depend on $g;$ this is the content of the following 
theorem, due to Galligo.
\bt
Let $I\subseteq SW$ be a homogeneous ideal. There exists a Zariski open subset
$U_I\subseteq\nR(W)$ such that $\forall g\in U_I\;\; \ini_g(I)$ is constant; 
it is called the {\em generic initial ideal} of $I,\;\gin(I).$\\
Furthermore, $\gin(I)$ is Borel-fixed.
\et
\pf \cite{g} 1.27.\qed\vspace{\baselineskip}\\
Let $I$ be a homogeneous ideal, then $I_{<d}$ (resp. $I_{\leq d}$) denotes the
ideal generated by the elements of $I$ of degree $<d$ (resp. $\leq d$); 
furthermore, we say that $\gin(I)$ has a {\it gap} in degree $\delta$ if 
$\gin(I)$ has no minimal generator of degree $\delta.$\\
The following result is due to Green.
\bp
If $\gin(I)$ has a gap in degree $\delta,$ then 
$\gin(I_{<\delta})=(\gin(I))_{<\delta}.$
\ep 
\pf \cite{ds} 2.12.\qed
\bc\label{crys}
If $\gin(I)$ has a gap in degree $\delta,$ then 
$\gin(I_{\leq\delta})=(\gin(I))_{\leq\delta}.$ 
\ec 
\pf Since $\gin(I)$ has a gap in degree $\delta,$ then 
$(\gin(I))_{\leq\delta}=(\gin(I))_{<\delta},$ because there is no generator in
degree $\delta,$ but also $I_{<\delta}=I_{\leq\delta},$ because a minimal 
Gr\"obner basis of $I$ (cf. \cite{clo} ch. 2) contains no polynomial of degree 
$\delta$ either. It follows that 
$\gin(I_{\leq\delta})=\gin(I_{<\delta})=(\gin(I))_{<\delta}=
(\gin(I))_{\leq\delta}.$
\qed\vspace{\baselineskip}\\
For ease of reference, we collect here a few facts that will be needed 
later on.
\bp{\normalshape  (\cite{g} 2.30)}\label{sat}
The following are equivalent:\\
(i) $I$ is saturated;\\
(ii) $\gin(I)$ is saturated;\\
(iii) no generator of $\gin(I)$ contains the last variable $x_n.$
\ep 
\bp{\normalshape (\cite{g} 2.14)}\label{slice}
For a general linear form $h,$ the equality 
$$\gin((I:h^p)|_h)=(\gin(I):x_n^p)|_{x_n}$$ 
holds for all $p\geq0.$
\ep 
\bt{\normalshape(\cite{s} theorem I.6.8)}\label{shaf}
Let $f:X\rightarrow B$ be a surjective morphism of projective varieties.
If $B$ is irreducible and all the fibers $f^{-1}(b)$ are irreducible and of 
the same dimension, then $X$ is irreducible.
\et 


\section{Connectedness}
If $W'\subseteq W$ is a linear subspace, the natural projection 
$W\rightarrow\frac{W}{W'}$ induces a map on the corresponding symmetric 
algebras $f\in\nS W\rightarrow f|_{W'}\in\nS\left(\frac{W}{W'}\right);$ 
especially, for any linear form 
$h\in W,\; f|_{h}\in\nS\left(\frac{W}{\langle h\rangle}\right).$\\
Note that $\frac{\kk\langle x_0,...,x_n\rangle}{\langle x_n\rangle}\simeq
\kk\langle x_0,...,x_{n-1}\rangle$ canonically, hence, for any 
$f\in\kk[x_0,...,x_n],$ we think of $f|_{x_n}$ as a polynomial in 
$\kk[x_0,...,x_{n-1}],$ and similarly, if ${\frak i}\subseteq\kk[x_0,...,x_n]$
is a (homogeneous) ideal, then ${\frak i}|_{x_n}$ is a (homogeneous) ideal of
$\kk[x_0,...,x_{n-1}].$\\ 
Let $S=V(W')\subseteq{\Bbb P}(W^*)$ be the subspace determined by $W',$ 
and let $Z\subseteq{\Bbb P}(W^*)$ be a projective variety, with ideal 
$I=I_Z,$ then, set-theoretically, $V(I|_{W'})=Z\cap S;$ especially, if $H=V(h)$
is a hyperplane, then $V(I|_h)=Z\cap H$ is the hyperplane section of $Z.$\\
Note that, for arbitrary $h,\; I|_h$ is a proper subset of $I_{Z\cap H};$
however, for general $h,$ their saturations coincide, i.e. for a general linear
form $h\in W,\; (I|_h)^{\isat}=I_{Z\cap H}$ (recall that $I_{Z\cap H}$ is 
prime, hence saturated).\\
Let ${\frak b}$ be a Borel-fixed monomial ideal in $\kk[x_0,...,x_n],$ then, 
for all ${\bf \tilde{p}}=(p_2,p_3,...,p_n)\in{\Bbb N}^{n-1}$ and 
${\bf \tilde{x}}^{{\bf \tilde{p}}}=x_2^{p_2}x_3^{p_3}\dots x_n^{p_n},$
the ideal $({\frak b}:{\bf \tilde{x}}^{{\bf \tilde{p}}})\cap \kk[x_0,x_1]$ 
is Borel-fixed too, hence it has a minimal system of generators of type
$x_0^s, x_0^{s-1}x_1^{\lambda_{s-1}},..., 
x_0x_1^{\lambda_1}, x_1^{\lambda_0},$ with $s=s({\bf\tilde{p}})$ and 
$\lambda_i=\lambda_i({\bf \tilde{p}}),$ for $i=0,1,...,s-1.$
\bd
The {\em monomial invariants} of ${\frak b}$ are the integers 
$\lambda_i({\bf \tilde{p}}),$ for all $i$ and ${\bf \tilde{p}}.$
\ed 
\brem
{\em
(i) By noetherianity, there are finitely many monomial invariants.\\
(ii) Borel-fixedness implies that 
$\lambda_i({\bf \tilde{p}})\geq\lambda_{i+1}({\bf \tilde{p}})+1,$ for all
${\bf \tilde{p}}.$\\
(iii) Note that $s({\bf 0})=$ smallest degree of a polynomial of ${\frak b};$
furthermore, if $p_j\leq q_j$ for all $j=2,...,n,$ then 
$s({\bf\tilde{p}})\geq s({\bf\tilde{q}}).$
}
\erem
Let $Z\subseteq {\Bbb P}^n$ be a projective variety and let $I_Z$ be its 
ideal, then no generator of $\gin(I_Z)$ contains the last variable $x_n$,
so its monomial invariants depend only on 
${\bf\hat{p}}:=(p_2,p_3,...,p_{n-1})\in{\Bbb N}^{n-2},$ i.e. they are 
independent of $p_n.$
\bd
The {\em monomial invariants} of $Z$ are the monomial invariants of 
$\gin(I_Z),$ i.e. the integers $\lambda_i({\bf \hat{p}}),$ for all 
${\bf\hat{p}}\in{\Bbb N}^{n-2},$ and $i=0,1,...,s({\bf\hat{p}})-1.$ 
\ed 
\brem
{\em
(i) For $p_{n-1}$ big enough, $\lambda_i({\bf \hat{p}})$ are the 
monomial invariants of the generic hyperplane section of $Z.$ Especially, if 
$Z$ has codimension two and all $p_j$ are big enough, $j=2,...,n-1,$ then 
$\lambda_i({\bf \hat{p}})$ are independent of ${\bf \hat{p}}$ and are the 
monomial invariants of the section $\Gamma=Z\cap \Pi$ of $Z$ with a general 
2--plane $\Pi.$\\
(ii) Note that $s({\bf 0})=s_Z=$ smallest degree of a polynomial vanishing on 
$Z$ and $s({\bf \hat{p}})=s_\Gamma=$ smallest degree of a polynomial vanishing 
on $\Gamma,$ when $\hat{p}_j\gg0$ for all $j.$  
}
\erem 
\bt\label{conn}
Let $Z$ be an integral nondegenerate projective variety of codimension two, 
and let $\Gamma$ be its generic section with a linear subspace of 
dimension two. 
If $s_Z=s_\Gamma,$ then the monomial invariants of $Z$ are {\em connected} 
i.e. they satisfy the following inequality:
$$
\lambda_{i+1}({\bf \hat{p}})+2\geq
\lambda_i({\bf \hat{p}})\geq
\lambda_{i+1}({\bf \hat{p}})+1.
$$
for all ${\bf\hat{p}}$ and $i=0,1,...,s({\bf\hat{p}})-2.$ 
\et
\pf
As already remarked, the second inequality is a consequence of Borel-fixedness,
so one needs to prove only the first inequality.
To avoid too cumbersome notations, we restrict ourselves to the case of a 
surface $\Sigma$ in ${\Bbb P}^4,$ the general case being completely similar.\\
As a preliminary step, recall that, in general, 
$s_\Sigma=s(0,0)\geq s(p,q)\geq s_\Gamma$ for all $p,q,$ hence the hypothesis
$s_\Sigma=s_\Gamma$ implies that all $s(p,q)$ are equal; denote this common 
value by $s.$\\ 
Let $I:=I_\Sigma$ be the ideal of $\Sigma$ and assume that the monomial 
invariants are not connected, i.e. there are integers 
$j, {\bar{p}}, {\bar{q}},$ with $0\leq j<s({\bar{p}},{\bar{q}})-1,$ such that 
$$\lambda_j({\bar{p}},{\bar{q}})>\lambda_{j+1}({\bar{p}},{\bar{q}})+2.$$
Set $\delta:=j+\lambda_{j+1}({\bar{p}},{\bar{q}})+2$ and define, for general 
linear forms $h, l, m,$  
$$J=J(h,l,m):=((I|_h:l^{\bar{q}})|_l:m^{\bar{p}})|_m.$$ 
{\bf Step 1}\; 
{\em $\gin(J)$ has a gap in degree $\delta.$}
\vspace{\baselineskip}\\
A repeated application of proposition \ref{slice}, shows that 
$$\gin(J)=(((\gin(I)|_{x_4}:x_3^{\bar{q}})|_{x_3}:x_2^{\bar{p}})|_{x_2}.$$
It follows that $\gin(J)$ is a Borel--fixed monomial ideal in 
$\kk[x_0,x_1]$ and has invariants 
$\lambda_i=\lambda_i({\bar{p}},{\bar{q}}),\, i=0,1,...,s-1;$ thus 
$\gin(J)=(x_0^s,x_0^{s-1}x_1^{\lambda_{s-1}},...,x_0^{j+1}x_1^{\lambda_{j+1}},
x_0^jx_1^{\lambda_j},...,x_1^{\lambda_0}),$ 
where the degrees of the generators are in increasing order; especially, since
$\deg (x_0^{j+1}x_1^{\lambda_{j+1}})=j+1+\lambda_{j+1}=\delta-1$ and 
$\deg (x_0^jx_1^{\lambda_j})=j+\lambda_j>j+\lambda_{j+1}+2=\delta,$ then 
$\gin(J)$ has a gap in degree $\delta.$ 
\vspace{\baselineskip}\\
{\bf Step 2}\; 
{\em There exists a homogeneous polynomial $F=F_{h,l,m}$ of degree $j+1$ 
such that 
$$J_{\leq\delta}\subseteq(F).$$
}
Since $\gin(J)$ has a gap in degree $\delta,$ then, by corollary \ref{crys}
$$\gin(J_{\leq\delta})=\gin(J)_{\leq\delta}=
(x_0^s,x_0^{s-1}x_1^{\lambda_{s-1}},...,x_0^{j+1}x_1^{\lambda_{j+1}})
\subseteq(x_0^{j+1}).$$
This relation shows that the Hilbert function $\gin(J_{\leq\delta}), \;P_{\gin(J_{\leq\delta})}(t)=j+1,$ for all $t\gg0;$ on the other hand, 
$J_{\leq\delta}$ and $\gin(J_{\leq\delta})$ have the same Hilbert 
function, hence $P_{J_{\leq\delta}}(t)=j+1$ definitively, so the variety 
$V(J_{\leq\delta})$ is a group of $j+1$ points in ${\Bbb P}^1;$ it 
follows that there is a homogeneous polynomial $F=F_{h,l,m}$ of degree $j+1$ 
such that $J_{\leq\delta}\subseteq(F).$ Note that,
as $j+1$ is the maximal degree of such a polynomial, $F$ is also the 
greatest common divisor (gcd) of the generators of $J_{\leq\delta}.$ 
\vspace{\baselineskip}\\
{\bf Step 3}\; 
{\em $F_{h,l,m}$ form an algebraic family.}
\vspace{\baselineskip}\\
It means that there exists a polynomial 
$F({\bf x},\xi,\eta,\theta),$ separately homogeneous in all 
groups of variables, of degree $j+1$ with respect to ${\bf x},$ such that, 
for general 
${\bf a},{\bf b},{\bf c}\in\kk^5,\,F({\bf x},{\bf a},{\bf b},{\bf c})$ 
restricted to 
$h={\bf a}\cdot{\bf x}, l={\bf b}\cdot{\bf x}, m={\bf c}\cdot{\bf x},$ is 
$F_{h,l,m};$ here ${\bf a}\cdot{\bf x}:=\sum_{i=0}^4a_ix_i,$ and similar 
meaning have ${\bf b}\cdot{\bf x}$ and ${\bf c}\cdot{\bf x}.$ This can be seen 
using the following argument.\\
Let $I=(\Phi_1,...,\Phi_r), \, \Phi_i\in\kk[x_0,...x_4],$ set 
$H=\xi\cdot{\bf x}, L=\eta\cdot{\bf x}, M=\theta\cdot{\bf x},$
and define $I|_H=(\Phi_1|_H,...,\Phi_r|_H),$ where 
$\Phi_i|_H=
\Phi_i(x_0,...,x_3,-\frac{\xi_0}{\xi_4}x_0...-\frac{\xi_3}{\xi_4}x_3)\in
\kk(\xi)[x_0,...,x_3]$ (of course, the substitution 
$x_4=-\frac{\xi_0}{\xi_4}x_0...-\frac{\xi_3}{\xi_4}x_3$ comes from 
$\xi\cdot{\bf x}=0$); with a similar substitution define the image of  
$L$ in $\kk(\xi,\eta)[x_0,...,x_3],$ image that we still denote by $L.$
Since the generators of the ideal quotient $(I_H:L^{\bar{q}})$ are 
algorithmically obtainable from 
the $\Phi_i$'s via rational operations in $\kk(\xi,\eta)[x_0,...,x_3]$ 
(cf. \cite{clo} theorem 4.4.1 and remark afterwards), we can therefore compute 
$(I_H:L^{\bar{q}})=(\Psi_1,...,\Psi_t),$
where $\Psi_i\in\kk(\xi,\eta)[x_0,...,x_3].$ Iterating the argument, we have 
${\frak J}:=((I|_H:L^{\bar{q}})|_L:M^{\bar{p}})|_M=(\Omega_1,...,\Omega_u),$ 
with $\Omega_i\in\kk(\xi,\eta,\theta)[x_0,x_1]$ and 
${\frak J}_{\leq\delta}=(\Omega_1,...,\Omega_v),\;v\leq u,$ taking only the 
generators of degree $\leq\delta.$ Now we compute 
${\rm gcd}(\Omega_1,...,\Omega_v)$ in $\kk(\xi,\eta,\theta)[x_0,x_1]$ using 
the euclidean algorithm; clear the denominators of 
${\rm gcd}(\Omega_1,...,\Omega_v)$ and the result is 
$F({\bf x},\xi,\eta,\theta).$ 
Summing up, $F({\bf x},\xi,\eta,\theta)$ is obtained from the 
$\Phi_i$ via an algorithm $\cal{A}$ that uses only rational operations; 
specializing $\xi\rightsquigarrow{\bf a}, \eta\rightsquigarrow{\bf b}, 
{\bf \theta}\rightsquigarrow{\bf c},$ where ${\bf a}, {\bf b}, {\bf c}$ 
are such 
that none of the denominators that appear in the alghorithm $\cal{A}$ vanishes 
(hence ${\bf a}, {\bf b}, {\bf c}$ are general in $\kk^5$), we get as a result 
$F({\bf x},{\bf a},{\bf b},{\bf c}),$ which, when restricted to $h, l, m,$ is 
the gcd of the generators of $J_{\leq\delta},$ hence 
$F({\bf x},{\bf a},{\bf b},{\bf c})|_{\langle h,l,m\rangle}=F_{h,l,m}.$
\vspace{\baselineskip}\\
Let $G\in I$ be a polynomial of minimal degree $s,$ then 
$G|_{\langle h,l,m\rangle}\in J_{\leq\delta},$ because $s\leq\delta,$ hence 
$F_{h,l,m}$ divides $G|_{\langle h,l,m\rangle},$ because 
$J_{\leq\delta}\subseteq(F_{h,l,m}),$ and furthermore it is a 
proper factor, for $\deg F_{h,l,m}=j+1<s=\deg G|_{\langle h,l,m\rangle}.$\\
Note that we use here the hypothesis $s_\Sigma=s_\Gamma;$ indeed, in general, 
$s({\bar{p}},{\bar{q}})<\delta,$ but, without the hypothesis 
$s_\Sigma=s_\Gamma,$ it could happen that $s(0,0)>\delta,$ so 
$G|_{\langle h,l,m\rangle}\not\in J_{\leq\delta}.$
\vspace{\baselineskip}\\
{\bf Step 4}\; 
{\em The set
$$X:=\{ (P,h,l,m)\in{\Bbb P}^4\times({\Bbb P}^{4*})^3 | P\in V(G,h,l,m) \}$$
is an irreducible variety.}
\vspace{\baselineskip}\\
Consider the canonical projections of $X$ on $V(G)$ and 
$({\Bbb P}^{4*})^3,$
\begin{equation}
\begin{CD}
X @>\pi_1>> V(G) \\
@V\pi_2VV \\
({\Bbb P}^{4*})^3
\end{CD}\nonumber
\end{equation}
The fiber of $\pi_1$ over a point $P\in V(G)$ consist of all 
$(h,l,m)\in({\Bbb P}^{4*})^3$ such that $P\in V(h,l,m),$ i.e. 
$\pi_1^{-1}(P)=\Pi_P\times\Pi_P\times\Pi_P,$ where 
$\Pi_P:=\{h\in{\Bbb P}^{4*} | P\in h\}$ is a hyperplane of ${\Bbb P}^{4*}$, 
hence all fibers of $\pi_1$ are irreducible and of the same dimension; $V(G)$ 
is irreducible too, because $G$ is a polynomial of minimal degree vanishing on 
an irreducible variety $\Sigma;$ thus $X$ is irreducible by theorem \ref{shaf}.
Note that $\dim X=\dim({\Bbb P}^{4*})^3=12,$ because the generic fiber of 
$\pi_2$ is zero dimensional.
\vspace{\baselineskip}\\
To conclude the proof of the theorem, let 
$Y:=\{ (P,h,l,m)\in{\Bbb P}^4\times({\Bbb P}^{4*})^3 | 
P\in V(h,l,m),F(P,h,l,m)=0 \}.$ Considering again the canonical projections
\begin{equation}
\begin{CD}
X\cap Y @>\pi_1>> {\Bbb P}^4 \\
@V\pi_2VV \\
({\Bbb P}^{4*})^3
\end{CD}\nonumber
\end{equation}
we notice that the general fiber $\pi_2^{-1}(h,l,m)$ is the zero set of 
$F({\bf x},{\bf a},{\bf b},{\bf c})=F_{h,l,m},$ i.e. $j+1$ points on the line 
$V(h,l,m)$ that belong to $V(G)$ too; thus 
$\dim X\cap Y=\dim({\Bbb P}^{4*})^3=12,$
and furthermore, since those points are only $j+1$ out of the $s$ points of
$V(G,h,l,m), \; X\cap Y$ is a proper subvariety of $X$ of the same dimension;
this of course contradicts the irreducibility of $X.$\\
This contradiction proves the theorem.\qed
\brem
{\em
To get a contradiction in the previous proof, one only needs that 
$\delta\geq s(0,0).$ This observation makes possible to push the 
statement of theorem \ref{conn} a little further, as follows. (For sake of 
simplicity, we only consider the case of a curve $C\subseteq{\Bbb P}^3.$)\\
(i) Even dropping the hypothesis $s_C=s_\Gamma,$ the invariants at level zero, 
i.e. the invariants $\lambda_0(p),$ are connected, because, at this level 
$\delta\geq s(0);$ the invariants at level one, $\lambda_1(p),$ are connected 
too, because $s(1)\geq s(0)-1$ by Borel-fixedness, so 
$\delta\geq s(1)+1\geq s(0);$ also, the invariants are connected at level 
$p\gg0,$ as already noticed in the introduction.\\ 
(ii) The monomial invariants of $C$ are likewise the monomial invariants of 
$I_C|_h:$ this follows from the fact that $\gin(I_C)$ has no generator 
containing $x_3$ (proposition \ref{sat}). 
Since $I_\Gamma=(I_C|_h)^\isat,$ for all 
$k\geq0,\, I_C|_h\subseteq (I_C|_h:{\frak m}^k)\subseteq I_\Gamma,$ where 
${\frak m}=(x_0,x_1,x_2)$ is the irrelevant maximal ideal of 
$\kk[x_0,x_1,x_2],$ and the 
hypothesis $s_C=s_\Gamma$ implies that a similar condition is satisfied also 
by $(I_C|_h:{\frak m}^k),$ and for such an ideal the proof of 
theorem \ref{conn} carries through, so we can conclude that the monomial 
invariants of $(I_C|_h:{\frak m}^k)$ are connected as well. 
}
\erem



\section*{}
e--mail: 
{\tt \{alzati, tortora\}{\char`@}mat.unimi.it}

\end{document}